\input amstex
\documentstyle{amsppt}

  \loadmsbm
  \loadbold

 \magnification=\magstep1
        \pagewidth{13cm}
        \pageheight{20cm}

\define\dF{{\Bbb F}}

\define\dQ{{\Bbb Q}}

\define\dZ{{\Bbb Z}}

\define\eq{{\frak q}}

\define\ep{{\frak p}}
\define\eP{{\frak P}}

\define\tw#1{\,^{#1}\!}

\font\pfeile = cmsy10 scaled 1440
\newfam\pfeilfam
\textfont\pfeilfam=\pfeile
                   \scriptfont\pfeilfam=\pfeile
                                      \scriptscriptfont\pfeilfam=\pfeile

\mathchardef\swpfeil="2D2E

\def\\{\let\stoken= } \\
\long\def\unexpandedwrite#1#2{\def\finwrite{\write#1}%
{\aftergroup\finwrite\aftergroup{\sanitize#2\endsanity}}}
\def\sanitize{\futurelet\next\sanswitch}
\def\sanswitch{\ifx\next\endsanity
\else\ifcat\noexpand\next\stoken\aftergroup\space\let\next=\eat
\else\ifcat\noexpand\next\bgroup\aftergroup{\let\next=\eat
\else\ifcat\noexpand\next\egroup\aftergroup}\let\next=\eat
\else\let\next=\copytoken\fi\fi\fi\fi \next}
\def\eat{\afterassignment\sanitize \let\next= }
\long\def\copytoken#1{\ifcat\noexpand#1\relax\aftergroup\noexpand
\else\ifcat\noexpand#1\noexpand~\aftergroup\noexpand\fi\fi
\aftergroup#1\sanitize}
\def\endsanity\endsanity{}

\expandafter\ifx\csname pre numero.tex at\endcsname\relax \else  \fi
\expandafter\chardef\csname pre numero.tex at\endcsname=\the\catcode`\@
\catcode`\@=11
\def\tokenitize{\futurelet\next\tokenswitch}
\def\tokenswitch{\ifx\next\endtokenity
\else\ifcat\noexpand\next\stoken\aftergroup\space\let\next=\eatt
\else\let\next=\copytok\fi\fi\next}
\def\eatt{\afterassignment\tokenitize \let\next= }
\long\def\copytok#1{\aftergroup\string\aftergroup#1\tokenitize}
\def\endtokenity{}
\newcount \secno
\newcount \Refno
\def\phantomlabel#1{{\aftergroup\expandafter\aftergroup\xdef%
\aftergroup\csname\tokenitize ref#1\endtokenity\aftergroup\endcsname}%
{{\the\secno.\the\Refno}}%
\global\advance\Refno by 1\relax}
\def\label#1{{\the\secno.\the\Refno}\phantomlabel{#1}%
\unexpandedwrite\labx{\phantomlabel{#1}}}
\def\Ref#1{{\aftergroup\expandafter\aftergroup\ifx
\aftergroup\csname\tokenitize ref#1\endtokenity}\endcsname\relax
??R\'ef\'erence {\tokenitize#1\endtokenity} non d\'efinie??%
\message{Reference #1 non definie}%
\else{\aftergroup\csname\tokenitize ref#1\endtokenity}\endcsname\fi
\ifnextsecn@\else\message{ATTENTION secno n'a pas ete initialise}\fi
}
\def\phantomsection{\global\advance\secno by 1 \global\Refno = 1\relax}
\def\secnum{{\phantomsection\the\secno\unexpandedwrite\labx{\phantomsection}}}
\newif\ifnextsecn@
\def\nextsecno#1{\nextsecn@true\global\secno=#1\global\advance\secno by -1%
\immediate\write\labx{\secno=\the\secno}\message{\the\secno}
}
\openin1 \jobname.lab \ifeof1 \message
                                {No file \jobname.lab.}\else\closein1\relax\input \jobname.lab \fi
\newwrite\labx
\immediate\openout\labx=\jobname.lab 
\catcode`\@=\csname pre numero.tex at\endcsname
     \nextsecno 1
\scrollmode \NoBlackBoxes \magnification=\magstep1
        \pagewidth{13cm}
        \pageheight{20cm}\topmatter
\define\AT{1}
\define\KS{2}
\define\matzat{3}
\define\saltman{4}
\define\VdW{5}
\define\weil{6}
\title
Irreducible polynomials which are locally reducible everywhere \endtitle
\author
Robert Guralnick, Murray M. Schacher and Jack Sonn \endauthor
\affil University of Southern California, Los Angeles,
University of California at Los Angeles, and
Technion--Israel Institute of Technology,
Haifa, Israel \endaffil

\address
R.~G.~:Department of Mathematics, University of Southern California, Los Angeles, CA
\endaddress \email guralnic\@usc.edu \endemail
\address
M.~M.~S.~:Department of Mathematics, University of California at Los Angeles, Los Angeles,
CA \endaddress \email mms\@math.ucla.edu \endemail
\address
J.~S.~:Department of Mathematics, Technion, 32000 Haifa, Israel \endaddress \email
 sonn\@math.technion.ac.il
\endemail

\thanks
The first author was partially supported by NSF Grant DMS 0140578. The research of the
third author was supported by Technion V.P.R. Fund--S. and N. Grand Research Fund
\endthanks \subjclass 11R52, 11S25, 12F05, 12G05, 13A20
 \endsubjclass
\abstract For any positive integer $n$, there exist polynomials $f(x)\in
\dZ[x]$ of degree $n$ which are irreducible over $\dQ$ and reducible 
over $\dQ_p$ for all
primes $p$, if and only if $n$ is composite.  In fact, this result holds
over arbitrary global fields. \endabstract \endtopmatter
\document

\head
  \secnum.  Introduction
\endhead
\smallskip
Hilbert  gave examples of irreducible polynomials $f(x)\in \dZ[x]$ 
of degree $4$
which
are reducible mod $p$ for all primes $p$, namely $x^4+2ax^2+b^2$. Note
that this polynomial is irreducible over $\dQ(a,b)$ hence (by Hilbert's
irreducibility theorem) is irreducible over $\dQ$ for infinitely many
specializations of $a$ and $b$ into $\dQ$.  The underlying reason for
this phenomenon from the Galois theoretic point of view is that the
Galois group of $x^4+2ax^2+b^2$ over $\dQ(a,b)$ is Klein's four group.
Therefore for any $p$ not dividing the discriminant of $f$, the
decomposition group is a cyclic group of order at most 2, so $f$ is
reducible mod $p$.  (Note that for $p$ dividing the discriminant of $f$,
$f$ is reducible mod $p$ as well.) The phenomenon is thus forced by the
structure of the Galois group. This also explains why there can be no
such examples of polynomials of prime degree. Indeed, suppose $f(x)\in
\dZ[x]$ has prime degree $\ell$ and is irreducible in $\dZ[x]$. Then its
Galois group has an element of order $\ell$, so by Chebotarev's density
theorem there exists $p$ such that the splitting field of $f$ over
$\dF_p$ has Galois group $C_{\ell}$, the cyclic group of order $\ell$,
hence $f$ must be irreducible over $\dF_p$. We will show that the degree
of $f$ being prime is the only obstacle, namely that for any composite
$n$, there exist irreducible $f(x)\in \dZ[x]$ of degree $n$ which are reducible mod
$p$ for all $p$. In fact, there is an irreducible $f(t,x)\in\dZ[t,x]$
of degree $n$
such that $f(t_0,x)$ is reducible mod $p$ for all specializations
$t=t_0$ in $\dZ$ and all $p$.  We will also prove that for any composite
$n$, there exist irreducible $f(x)\in \dQ[x]$ of degree $n$ which are reducible over
$\dQ_p$ for all $p$, and that this result generalizes to arbitrary
global fields. Note that Hilbert's example does not satisfy this last
condition for all $a,b$, e.g. $x^4+1$ is irreducible over $\dQ_2$.
\smallskip It is worthwhile pointing out here that a random polynomial $f(x)\in\dZ$ of composite
degree $n$ is not reducible mod $p$ for all $p$, as its Galois group
over $\dQ$ is $S_n$ \cite\VdW, and since $S_n$ contains an $n$-cycle,
Chebotarev's density theorem implies that there are infinitely many
primes $p$ for which $f(x)$ is irreducible mod $p$.  
\bigskip
\head \secnum. mod $p$ reducibility \endhead

Let $f(x)\in \dZ[x]$ be monic irreducible of degree $n$ with Galois
group $G$ over $\dQ$. As in the discussion above, we see that if $f(x)$
is irreducible mod $p$, then $p$ does not divide the discriminant of
$f(x)$, and $G$ must contain an element of order $n$, since the
decomposition group is cyclic of order $n$. Thus if $G$ has no element
of order $n$, then $f$ must be reducible mod $p$ for all $p$. On the
other hand, since $f$ is irreducible over $\dZ$,  $G$ has a subgroup $H$
of index $n$. More generally, if $K/\dQ$ is a finite Galois extension
with Galois group $G$, and if $G$ has a subgroup $H$ of index $n$ but no
element of order $n$, then the fixed field, say $\dQ(\theta)$, of $H$
has degree $n$ and the minimal polynomial $f$ of $\theta$ over $\dQ$ has
degree $n$ with splitting field $L\subseteq K$, Galois group $\bar G =
G(L/\dQ)\cong G/core(H)$ where $core(H)$ is the intersection of the
conjugates of $H$.  Furthermore, $\bar G$ has a subgroup $\bar H$ of
index $n$ but has no element of order $n$, so $f$ is reducible mod $p$
for every $p$. In summary:
\smallskip
\proclaim {Lemma \label{basic}} Let $G$ be a finite group and $n$ a positive integer such
that

{\rm 1)}\ \ \ \ \ $G$ is realizable as a Galois group $G(K/\dQ$,

{\rm 2)} \ \ \ \ \ $G$ has a subgroup of index $n$,  but $G$ has no
element  of order $n$.

Then there is an irreducible polynomial $f(x)\in \dZ[x]$ of degree $n$
which is reducible mod $p$ for all primes $p$ (with splitting field
contained in $K$). \endproclaim

\bigskip
In connection with Hilbert's example, it follows from the preceding
considerations that a polynomial of degree four has the property of
being irreducible over $\dZ$ and reducible mod $p$ for all $p$ if and
only if its Galois group over $\dQ$ is either Klein's four group $V_4$
or the alternating group $A_4$.

 \proclaim {Theorem \label{modp}} For any composite positive
integer $n$, there exist irreducible polynomials $f(x)\in \dZ[x]$ of
degree $n$ which are reducible mod $p$ for all primes $p$. \endproclaim
\smallskip
\demo{Proof} By the lemma, it suffices to find $G$ satisfying conditions 1) and 2).
\smallskip
\it Case 1. \rm $n$ squarefree.
Write $n=qm$ with $q$ prime, so $(m,q)=1$.
Let $t$ be the order of $q$ mod $m$ and set $G=C_m\ltimes V$ (semidirect
product), where $V$ denotes the additive group of $\dF_{q^t}$ with $C_m$
acting by multiplication by the $m$th roots of unity  $ \mu_m
\subset\dF_{q^t}$, where we identify $C_m$ with $\mu_m$.  The action of
$C_m$ on $V$ is faithful and $0$ is the only fixed point.  Let $H$ be a
hyperplane in $V$, i.e. $[V:H]=q$.  Then  $[G:H]=qm=n$.  Furthermore, if
$g\in G$ had order $n$, then $g^m$ ($\in V$) would have order $q$ and be
fixed by $g$ hence by $C_m$, contradiction.


We therefore have found for every squarefree composite $n$, a solvable
group $G$ satisfying condition 2).  By Shafarevich's theorem (in fact by
an older theorem of Scholz), $G$ is realizable as a Galois group over
$\dQ$.
\smallskip
\it Case 2. \rm $n$ not squarefree.  Assume first that $n=q^2$, $q$ prime.  Then
$G_1:=C_q\times C_q$ satisfies 2) with $H$ the trivial group.  For arbitrary $n=q^2m$, take
$G:=G_1\times C_m$ which again has the trivial subgroup of index $n$ and no element of
order $n$. \qed \enddemo
\medskip
We now remark that since the groups $G$ appearing in the proof can be
realized regularly over the rational function field $\dQ(t)$ (see e.g.
\cite {\matzat, p. 275}), we obtain irreducible polynomials $f(t,x)\in
\dZ[t,x]$ of degree $n$ which by Hilbert's irreducibility theorem have
infinitely many specializations of $t$ into $\dQ$ which are irreducible
with Galois group $G$, hence reducible mod $p$ for all $p$.

\head \secnum. $p$-adic reducibility \endhead

We now wish to prove Theorem~\Ref{modp} with reducibility mod $p$
replaced by reducibility over $\dQ_p$.  The preceding construction
actually yields irreducible polynomials $f(x)\in \dZ[x]$ which are
reducible over $\dQ_p$ for all primes which are unramified in the
splitting field of $f$, but may be irreducible over $\dQ_p$ for ramified
$p$.  The proof will be similar but more delicate.
\smallskip
Let $f(x)\in \dQ[x]$ be irreducible and let $p$ be a prime.  Let $K$ be the
splitting field of $f$ over $\dQ$, $G=G(K/\dQ)$, $\ep$ a prime of $K$ over
$p$, $D=D(\ep)$ the decomposition group.  Then $f$ is irreducible over
$\dQ_p$ if and only if $D(\ep)$ acts transitively on the roots of $f$ in
$K$.  Let $H$ be the subgroup of $G$ fixing a root of $f$.  The action of
$G$ on the roots of $f$ is equivalent to its action by multiplication from
the left on the left cosets of $H$ in $G$.  Thus $f$ is irreducible over
$\dQ_p$ if and only if $D(\ep)$ acts transitively on the left cosets of $H$
in $G$, i.e. the set product $DH$ is equal to $G$.  Suppose $K/\dQ$ were
tamely ramified, so that all the decomposition groups were metacyclic.  Then
in order to insure that $f$ is reducible over all $\dQ_p$, it would suffice
that (the set product) $MH\neq G$ for every metacyclic subgroup $M$ of $G$.
We therefore have the following lemma.
\medskip
\proclaim {Lemma \label{basic1}} Let $G$ be a finite group and $n$ a positive integer such
that

{\rm 3)} \ \ \ \ $G$ is realizable over $\dQ$ by a tamely ramified
extension $K/\dQ$,

{\rm 4)} \ \ \ \ $G$ has a subgroup $H$ of index $n$, the intersection of whose conjugates
is trivial,  such that the set product $MH\neq G$ for every metacyclic subgroup $M$ of $G$.

Then there exists an irreducible polynomial $f(x)\in \dQ[x]$ of degree
$n$ (with splitting field $K$) which is reducible over $\dQ_p$ for all
primes $p$. \endproclaim
\medskip
\proclaim {Theorem \label{p-adic}} For any composite positive integer $n$, there exist
irreducible polynomials $f(x)\in \dQ[x]$ of degree $n$ which are reducible over $\dQ_p$ for
all primes $p$. \endproclaim \smallskip

\demo{Proof} By Lemma~\Ref{basic1}, it suffices to find $G$ satisfying 3) and 4).

\it Case 1. \rm $n$ squarefree.  The same $G$ and $H$ work as in the proof of
Theorem~\Ref{modp}, provided the order $t$ of $q$ mod $m$ is greater than $1$ (which can be
ensured by taking $q$ to be the smallest prime dividing $n$).  Indeed, $V$ has
$\dF_q$-dimension $t$ and is irreducible as a $C_m$-module, so in this case the only
metacyclic subgroups of $G$ are cyclic of order $q$ or a divisor of $m$. Thus condition 4)
is satisfied.
Condition 3) can be proved using Saltman's results \cite{\saltman, Theorem 3.5} to verify
that $G$ has a generic Galois extension and then \cite{\saltman, Theorem 5.9} to see that
there exists a tame extension realizing $G$.
\medskip
\it Case 2.  \rm $n$ not squarefree.  As in the proof of Theorem~\Ref{modp}, we first
assume $n=q^2$, $q$ a prime.  We will construct a Galois extension $K/\dQ$ with Galois
group $C_q\times C_q$, with local degree 1 or $q$ at all primes $p$.  If $K=\dQ(\theta)$
with $f(x)$ the minimal polynomial of $\theta$ over $\dQ$, then $f$ has the desired
property. As we saw earlier, for any prime unramified in $K$ the decomposition group is
cyclic, so the local degree is 1 or $q$ just from the structure of $G$.  We need to
construct $K$ so that the local degree is $q$ at the ramified primes. For this we use an
idea from \cite\KS.

Let $\ell$ be a prime congruent to 1 mod $q$ and let $L_{\ell}\subseteq \dQ(\mu_{\ell})$
such that $G(L_{\ell}/\dQ)\cong C_q$.  We seek another prime $r$ such that

5) \ \ \ \ $r\equiv 1$ \ \ \ (mod $q$)

6) \ \ \ \ $r$ splits completely in $L_{\ell}$

7) \ \ \ \ $\ell$ splits completely in $L_r$

\smallskip
If $r$ satisfies these conditions, then clearly $K=L_{\ell}L_r$ will have the local degree
at most $q$ everywhere.

5) is equivalent to the condition that $r$ splits completely in $\dQ(\mu_q)$.  7) is
equivalent to the condition that the Frobenius automorphism of $\ell$ in
$G(\dQ(\mu_r)/\dQ)$ fixes $L_r$ pointwise, i.e. is a $q$th power in the cyclic group
$G(\dQ(\mu_r)/\dQ)$, which is equivalent to $\ell$ being a $q$th power mod $r$, which in
turn means that the polynomial $x^q-\ell$ has a root mod $r$.  Since $\dF_r$ contains the
$q$th roots of unity, this is equivalent to $x^q-\ell$ factoring into linear factors mod
$r$, which is equivalent to the condition that $r$ splits completely in
$\dQ(\mu_q,\root{q}\of{\ell})$.  It follows that conditions 5),6),7) together are
equivalent to the condition that $r$ splits completely in
$L_{\ell}(\mu_q,\root{q}\of{\ell})$.  By Chebotarev's density theorem, such an $r$ exists.
This completes the case $n=q^2$.

Now assume the general case $n=q^2m$, $q$ prime.  Let $K$ be as in the case $n=q^2$, and
let $L$ be any abelian extension of $\dQ$ of degree $m$ such that $K\cap L=\dQ$.  (For
example, choose a prime $p\equiv 1$ (mod $m$), $p\neq \ell,r$, and let $L\subseteq
\dQ(\mu_p)$ of degree $m$ over $\dQ$.)  $KL/\dQ$ has degree $n$, and the local degree at
any prime is at most $qm<n$. \qed \enddemo
\smallskip
\head \secnum. Global fields \endhead

Theorem \Ref{p-adic} generalizes to arbitrary global fields $F$.  If $n$
is prime and $f(x)\in F[x]$ is a separable monic irreducible polynomial,
then by Chebotarev's density theorem, which holds over any global field
\cite {\weil, p. 289}, there exist primes $\ep$ of $F$ such that $f$ is
irreducible over the completion $F_{\ep}$.  (Even if $n$ is equal to the
characteristic $p$ of $F$, and $f(x)$ is inseparable of degree $p$, i.e.
$f(x)=x^p-a$, then by \cite {\AT, Chapter 9, Theorem 1}, if $f(x)$ is
irreducible, then $f(x)$ is irreducible over $F_{\ep}$ for infinitely
many $\ep$.)

\proclaim {Theorem \label{global}}  For any composite positive integer $n$, and any global
field $F$,  there exist irreducible polynomials $f(x)\in F[x]$ of degree $n$ which are
reducible over $F_{\ep}$ for all primes $\ep$ of $F$. \endproclaim \smallskip

\demo {Proof} If $F$ is a number field, one can reduce the proof to the case
$F=\dQ$.   Let $f(x)\in \dQ[x]$ be irreducible of degree $n$ and reducible
over $\dQ_p$ for all $p$, and suppose its splitting field $K$ satisfies
$K\cap F=\dQ$.  Then $f(x)$ is irreducible over $F$ and is reducible over
$F_{\ep}$ for all primes $\ep$ of $F$.  It remains to observe that the proof
of Theorem~\Ref{p-adic} produces infinitely many $\dQ$-linearly disjoint
extensions $K$ with the desired properties.
\smallskip
We may therefore assume that $F$ a global function field of
characteristic $p$. \smallskip

Let $f(x)\in F[x]$ be monic irreducible of degree $n$, fix a root
$\alpha$ in a splitting field $K$, let $G=G(K/F)$, and let
$H=G(K/F(\alpha))$. Let $\ep$ be a prime of $F$, $\eP$ a prime of $K$
dividing $\ep$, $D=D(\eP)$ the decomposition group.  Then $f(x)$ is
reducible over $F_{\ep}$ if and only if $D$ does not act transitively on
the roots of $f(x)$, i.e. the set $DH$ is not equal to $G$.  We reduce
the proof to the case $n$ a product of two primes, not necessarily
distinct.
\smallskip
For any positive integer $m$, let $E/F$ be a Galois (e.g. cyclic)
extension of degree $m$ such that $K\cap E=F$, with $K$ as above. Assume
for all primes $\eP$ of $K$ that $D(\eP)H\neq G$.  Let $\hat K=KE$,
$\hat G=G(KE/F)\cong G\times G(E/F)$.  Identify $H$ with $H\times
\{1\}\subset \hat G$.  Then for any prime $\hat \eP$ of $\hat K$,
$D(\hat \eP)H\neq \hat G$.  This reduces the proof of
Theorem~\Ref{global} to the case $n$ a product of two primes, say $r$
and $s$.\smallskip We will further reduce the proof to the case $F$ is a
rational function field $\dF_q(t)$, by constructing the desired
$K/\dF_q(t)$ linearly disjoint over $\dF_q(t)$ from any $F$ given in
advance.  Accordingly we now assume $F=\dF_q(t)$.   As we will see, when
$n$ is prime to $p$, one can give a proof which is analogous to that of
Theorem~\Ref{p-adic}. \smallskip
  \it Case 1. \rm $r=s\neq p$. We use a
function field analogue of a construction in \cite \KS. Let $\ep=(h(t))$
be a (finite) prime of $F$ ($h(t)\in \dF_q[t]$) which splits completely
in the extension $F'$ of $F$ obtained by adjoining all $r$th roots of
all elements of $\dF_q$ (including the $r$th roots of unity).  Then
$Cl_F(\ep)$, the ray class group mod $\ep$, is cyclic (the class group
of $F$ is trivial), and has order divisible by $r$. There is a
corresponding ray class field extension $R^{\ep}/F$ which is geometric
(regular over $\dF_q$), $G(R^{\ep}/F)\cong Cl_F(\ep)$, and $\ep$ is the
only prime that ramifies in $R^{\ep}$.  Let $L^{\ep}$ be the (unique,
cyclic) subfield of $R^{\ep}$ of degree $r$ over $F$.  We seek $\ep,\eq$
such that $K:=L^{\ep}L^{\eq}$ has the desired property, which is that
$\ep$ splits completely in $L^{\eq}$ and $\eq$ splits completely in
$L^{\ep}$.  Let $\ep$ be as above, which exists e.g. by Chebotarev's
density theorem (we don't really need Chebotarev yet since the extension
is a constant extension, but presently we will need it). We seek $\eq$
which satisfies the same conditions as does $\ep$, and additionally,
$\ep$ splits completely in $L^{\eq}$ and $\eq$ splits completely in
$L^{\ep}$. The last condition is an additional Chebotarev condition on
$\eq$ which is compatible with the preceding ones, so can be satisfied
by Chebotarev's density theorem. The condition $\ep$ splits completely
in $L^{\eq}$ is equivalent to the condition that the Frobenius
$Frob(\ep)$ is an $r$th power in $G(R^{\eq}/F)$, which is equivalent to
$\ep=h(t)$ being an $r$th power in $Cl_F(\eq)$, which is equivalent to
$h(t)$ being an $r$th power in the multiplicative group of the residue
field $\dF_q[t]/\eq$, which is equivalent to $\eq$ splitting completely
in $F(\root{r}\of{h(t)})$. This is another Chebotarev condition on $\eq$
compatible with the preceding ones.  The construction yields infinitely
many linearly disjoint extensions $K$, completing the proof of this
case.
\smallskip
\it Case 2. \rm $r\neq s$, both $\neq p$. Here we can argue exactly as in
the case $n$ squarefree in the proof of Theorem~\Ref{p-adic}, by
constructing a tamely ramified $G$-extension  using Saltman's results as
quoted earlier.  This yields infinitely many linearly disjoint extensions.
\smallskip
\it Case 3. \rm $r=s=p$.  We seek a pair of Artin-Schreier extensions
$L,M$ defined by $x^p-x-g(t)$,  $x^p-x-h(t)$ respectively, where $g(t)$,
$h(t)\in F$, such that $K:=LM$ has local degree $1$ or $p$ everywhere.
Let $g(t)$ be a polynomial without constant term, i.e. $g(0)=0$.  Since
the derivative of $x^p-x-g(t)$ with respect to $x$ is $-1$, the  only
prime of $F$ that ramifies in $L$ is infinity.  Furthermore, the prime
corresponding to $t$ splits completely in $L$ since $x^p-x-g(t)$ has $p$
distinct roots mod $t$.  The automorphism $t\mapsto 1/t$ interchanges
$(t)$ and infinity.  Let $h(t)=g(1/t)$.  Then infinity splits completely
in $M$ and $(t)$ is the only ramified prime in $M$.  It follows that the
local degree of $LM$ is $1$ or $p$ everywhere, as desired.  We get
infinitely many linearly disjoint extensions by varying $g(t)$, e.g. by
taking $g(t)$ of the form $t^e$, with $e$ prime to $p$.  The genus of
the corresponding curve grows with $e$, or one can show that there are
infinitely many such distinct extensions even over the completion of $F$
at infinity.
\smallskip
\it Case 4. \rm $r=p\neq s$.
\smallskip
\ \ \ \ \ Subcase 4.1. \ \ \ $p\nmid s-1$.   We can use the group $G$ we
used earlier. Let $C_p$ act on an irreducible $\dF_s$-space $V$, so
$dim(V)>1$. Claim that $ G:= C_p \ltimes V$ cannot be a local Galois group
anywhere. Indeed, if it were, the extension would be ramified, so the
inertia group would be a nontrivial normal subgroup of $G$, hence equal to
$V$ or $G$. There is no wild ramification, since $G$ has no normal
$p$-subgroup. But then the inertia subgroup must be cyclic, contradicting
the fact that $G$ has no cyclic normal subgroup. This proves the claim.  It
follows that all decomposition groups are proper subgroups of $G$.  But if
$D$ is a proper subgroup of $G$, then $DH$ cannot equal $G$, as is easily
verified (see beginning of the previous section). It remains to show that
$G$ is realizable (infinitely often, linearly disjointly) over $F$. Realize
$C_p$ by an Artin-Schreier extension $L/F$ (there are infinitely many), and
consider the embedding problem with kernel $V$. It has a proper solution
since $F$ is Hilbertian \cite {\matzat, p. 275}.
\smallskip
 \ \ \ \ \ Subcase 4.2. \ \ \ $p|s-1$.
\smallskip
Here the roles of $p$ and $s$ need to be interchanged, and we have to worry
about wild ramification.  (The argument will also cover subcase 4.1.)  Let
$E/F$ be a cyclic regular extension of degree $s$.
Claim there exists an Artin-Schreier extension $L/E$ of degree $p$ such that
every prime $\eq$ of $E$ that ramifies in $L$ is split completely over $F$
and all its remaining $s-1$ conjugates over $F$ are unramified in $L$.
Indeed, there are infinitely many primes $\ep$ of $F$ that split completely
in $E$.  Let $\{\ep_i\}$  be a sequence of such primes and let $\eq_i$ be a
prime of $E$ dividing $\ep_i$ for each $i$.  Let $h=h_E$ be the class number
of $E$.  Then in the sequence
$\eq_1,\eq_1\eq_2,...,\eq_1\eq_2\cdots\eq_{h+1}$ of ideals of $E$, there
must be two which differ by a principal ideal $(f)$; i.e. there exist $i\leq
j\in \{1,2,...,h+1\}$ such that $\eq_i\cdots\eq_j=(f)$ with $f\in E^*$. Take
$L/E$ to be the Artin-Schreier extension defined by the equation
$x^p-x=1/f$.  Since the only poles of $1/f$ are $\eq_i,...,\eq_j$, the only
ramified primes in $L/E$ are also $\eq_i,...,\eq_j$.  Furthermore, since
these are simple poles, the extension $L/E$ is nontrivial (of degree $p$).
This proves the claim.
\smallskip Let $M/F$ be the Galois closure of $L/F$, and set $W=G(M/E)$.
By construction, the conjugates of $L$ over $F$ are linearly disjoint over $E$, so $W$ is
isomorphic to the group ring $\dF_pC_s$ as $C_s$-modules, identifying $G(E/F)$ with $C_s$.
There is an irreducible submodule $V$ of $W$ on which $C_s$ acts faithfully.  The
$\dF_p$-dimension of $V$ is necessarily greater than 1 since $s>p$. Let $K$ be the subfield
of $M$ corresponding to the complementary submodule to $V$ in $W$ (by Maschke's theorem).
Then $K/F$ is Galois with group $G=C_s\ltimes V$.
  Let $H$ be a subgroup of $V$ of index $p$. We show that $H$ has the
desired property, that $DH\neq G$ for any decomposition group $D$.  The argument is similar
to previous ones.  If $D$ is the decomposition group of an unramified prime, then $D$ is
cyclic.  Since $G$ has no cyclic subgroup of order $sp$, $D$ has order $p$ or $s$, hence
$DH$ cannot equal $G$.  At a tamely ramified prime, the inertia group must be of order $s$
and be normal in $D$.  So $D$ is either cyclic of order $s$, and $DH\neq G$, or $D$ is of
order $sp$ with a normal $p$-Sylow subgroup, which cannot happen because $D$ has no element
of order $sp$.  Finally, the decomposition group of a wildly ramified prime in $K/F$ must
be contained in $V$, since such a prime must divide one of $\eq_i,...\eq_j$, and $C_s$ does
not fix any of them.  Hence again $DH\neq G$.  Finally, there are infinitely many linearly
disjoint extensions $E/F$, hence infinitely many linearly disjoint extensions $K/F$. The
proof is complete. \qed \enddemo

 \enddocument